\documentclass[12pt,a4paper]{article}

\usepackage[utf8]{inputenc}
\usepackage[english]{babel}
\usepackage{amsthm}
\usepackage{amsfonts}
\usepackage{amssymb}
\usepackage{amsmath}
\usepackage[affil-it]{authblk}
\usepackage{mathrsfs}
\usepackage[usenames]{color}
\usepackage{graphicx}

\newtheorem{Theorem}{Theorem}

\newtheorem{Lemma}{Lemma}
\newtheorem*{LemmaNoCounter}{Lemma}
\newtheorem{Statement}{Statement}

\newenvironment{Def}{\par\noindent{\bf Definition. }}{\par}
\newenvironment{Note}{\par\noindent{\bf Note. }}{\par}

\DeclareMathOperator{\re}{Re}
\DeclareMathOperator{\im}{Im}

\DeclareMathOperator{\codim}{codim}

\DeclareMathOperator{\const}{const}
\DeclareMathOperator{\supp}{supp}
\DeclareMathOperator*{\esssup}{ess\,sup}

\newcommand{\CC}{\mathbb{C}}
\newcommand{\RR}{\mathbb{R}}
\newcommand{\NN}{\mathbb{N}}

\newcommand{\vint}[2]{\bigl|_{#1}^{#2}\bigr.}
\newcommand{\sint}[3]{\{#1\}_{#2}^{#3}}
\newcommand{\set}[2]{\bigl\{#1\,\bigl|\bigr.\,#2\bigr\}}

\makeatletter
\renewcommand{\section}{\@startsection{section}{1}%
{0pt}{3.5ex plus 1ex minus .2ex}%
{2.3ex plus.2ex}{\normalfont\large\bfseries}}
\makeatother

\title{Conditions for sectoriality and compactness of the resolvent for a non-self-adjoint Sturm--Liouville operator with singular distributional potential}
\author{Sergey N. Tumanov}
\affil{Moscow Center of Fundamental and Applied Mathematics. Lomonosov Moscow State
University,
Moscow, Russia}

\begin{document}

\maketitle

\begin{abstract}
The aim of this paper is to find necessary and sufficient conditions for sectoriality and compactness 
of the resolvent for Sturm--Liouville operators with complex-valued potentials of the class $q\in W_{2,loc}^{-1}(\RR_+)$ in terms of its generalized antiderivatives $s\in L_{2,loc}(\RR_+)$.
\end{abstract}
\frenchspacing

\section{Main results}

In the classical theory of the Sturm--Liouville operator defined by the differential expression
$$
l(y)=-y''+qy,
$$
where $q\in L_{1,loc}(\RR_+)$,
general results that do not require additional smoothness of the potential $q$ are often
formulated in terms of its antiderivative $s$, where $q=s'$. In particular, the Molchanov criterion for the compactness of the resolvent \cite{Molchanov} and its generalizations
\cite{Brinck,Lidsky,Ismagilov,TumMatSbor24} characterize the growth of the antiderivative on segments of constant length.

The nature of this phenomenon can be explained by the possibility of regularization of the differential equation $l(y)=f$ for $f\in L_{1,loc}(\RR_+)$: 
there is an equivalent system of equations, the coefficients of which are functions of the antiderivative:
$$
\begin{pmatrix}
y(x)\\
y^{[1]}(x)
\end{pmatrix}'=
\begin{pmatrix}
s(x) & 1 \\
-s^2(x) & -s(x)
\end{pmatrix}
\begin{pmatrix}
y(x)\\
y^{[1]}(x)
\end{pmatrix}-
\begin{pmatrix}
0\\
f(x)
\end{pmatrix}.
$$
The equivalent notation for the expression $l(y)$ is following:
\begin{equation}
\label{eqly}
l(y)=-(y^{[1]})'-sy^{[1]}-s^2y,\ \mbox{where } y^{[1]}=y'-sy.
\end{equation}

It is expected that the general spectral properties of the regularized problem are expressed in terms of $s$.

In their pioneering work \cite{SavchShkal} Savchuk and Shkalikov showed that regularization is an effective way of
generalizing the Sturm--Liouville operator for potentials $q\in W_{2,loc}^{-1}(\RR_+)$. Here $W_{2,loc}^{-1}(\RR_+)$ 
is the space of distributions with generalized antiderivatives $s\in L_{2,loc}(\RR_+)$.

The history of the issue and references are presented in \cite{SavchShkal, EGNT, mirzTMMO}, on the regularization method and its development
for high-order operators we recommend \cite{MirzShkal}.

Here we consider an operator with a complex-valued potential, and our aim is to obtain generalizations of the known classical
criteria for the compactness of the resolvent in terms of the antiderivative $s$ of the potential $q\in W_{2,loc}^{-1}(\RR_+)$. 
Throughout what follows $q=s'$, and the derivative is understood in the sense of distribution theory.

In a recent paper \cite{TumMatSbor24} we proposed a necessary condition for the compactness of resolvents for operators with arbitrary complex-valued potentials.
$q\in L_{1,loc}(\RR_+)$:
$$
\forall a>0\ \lim\limits_{x\to+\infty}\int\limits_x^{x+a}|q(\xi)|\,d\xi=+\infty.
$$

Of course, it is not relevant for $q\not\in L_{1,loc}(\RR_+)$, and here we make up for this deficiency by proposing Theorem \ref{thness} 
--- a stronger condition covering operators with arbitrary complex-valued $q\in W_{2,loc}^{-1}(\RR_+)$.

In \cite{TumMatSbor24} we found the
operator $L$ with essential spectrum, with potential $q\in L_{1,loc}(\RR_+)$, $\re q\ge0$ and $|q|\to+\infty$. The numerical range of $L$ is contained in the 
right half-plane, and the resolvent is bounded in the left half-plane.

It turned out that for a successful generalization of classical theorems of Molchanov type to the complex-valued case, one should require strict 
sectoriality for the minimal operator $L_0$ --- when the central angle of the sector with the numerical range of $L_0$ is strictly less than $\pi$.

For this reason, in the present paper we pay special attention to sectorial operators, for which we propose the corresponding theorems.

In the special case of a real-valued potential of the form:
\begin{equation}
\label{eqdeltasum}
q(x)=\sum_{x_k\in X}\alpha_k\delta(x-x_k),
\end{equation}
and more general case, when the antiderivative $s$ is a function of bounded
variation, the issues of semi-boundedness and compactness of the resolvent have been studied previously. For example, Brasche \cite{Brasche}, arguing 
in the language of quadratic forms,
proposed a sufficient condition for semi-boundedness of the quadratic forms:
$$
\int\limits_\RR|y'(x)|^2\,dx+\int\limits_\RR|y(x)|^2\,ds(x),\quad y\in W_2^1(\RR)\cap L_2(\RR,s),
$$
where $s$ is (neither necessarily finite nor necessarily positive) a Radon measure on $\RR$. For the form to be semi-bounded, it is sufficient:
\begin{equation}
\label{eqbrache}
\sup_{x>0}\bigl(s^{-}(x+1)-s^{-}(x)\bigr)\le C,
\end{equation}
where $s=s^{+}-s^{-}$ is a decomposition into non-negative measures on $\RR$.

A generalization of Molchanov's criterion for potentials of the form \eqref{eqdeltasum} was made in \cite{AlbKostMal}. 
Proposing a criterion for compactness of the resolvent for semi-bounded operators, the authors require the Brachet condition \eqref{eqbrache}.

The classical result of Brink \cite{Brinck} suggests a possible weakening of this requirement to:
$$
\sup_{x>0}\bigl(s(x+a)-s(x)\bigr)\ge -C\ \mbox{for all }0<a\le1
$$
In fact, the component $s^{+}$ can compensate for an arbitrarily large $s^{-}$, ensuring semi-boundedness of the operator.
This observation follows from theorem \ref{thbrink} below.

In \cite{KostenkoMalamud} the semi-boundedness of operators with potentials of the form \eqref{eqdeltasum} and discreteness of the spectrum were studied
using the theory of Jacobi matrices. Note that the proposed sufficient condition for semi-boundedness \cite[Cor.5.26]{KostenkoMalamud} follows from our theorem \ref{thbrink}.

The present work is substantially inspired by the results of Ismagilov \cite{Ismagilov} and Zelenko \cite{Zelenko},
and, of course, by the already mentioned work of Savchuk and Shkalikov \cite{SavchShkal}.
\bigskip

Let $s\in L_{2,loc}(\RR_+)$ be a complex-valued function.
For $y\in AC_{loc}(\RR_+)$, we define
$$
y^{[1]}(x)=y'(x)-s(x)y(x).
$$
For $y^{[1]}\in AC_{loc}(\RR_+)$ the following differential expression is correctly defined:
$$
l(y)=-(y^{[1]})'-sy^{[1]}-s^2y.
$$

We introduce lineals:
$$
\begin{aligned}
&\mathscr{D}=\set{y\in L_2(\RR_+)}{y, y^{[1]}\in AC_{loc}(\RR_+),\,l(y)\in L_2(\RR_+)},\\
&\mathscr{D}_0=\set{y\in\mathscr{D}}{y(0)=y^{[1]}(0)=0\,, \supp y\mbox{ --- is a compact}},
\end{aligned}
$$
and define operators in $L_2(\RR_+)$: the {\it minimal} $L_0$ on the domain $\mathscr{D}_0$ and the {\it maximal} $L_M$ on the domain $\mathscr{D}$,
where:
\begin{equation}
\label{eqL0minop}
L_0y=l(y),\, y\in \mathscr{D}_0,\quad L_My=l(y),\, y\in \mathscr{D}.
\end{equation}

Let us make a few reservations regarding the terminology used.

Below $s\in L_{2,loc}(\RR_+)$ will be called the {\it generalized antiderivative of the potential}, or simply the {\it antiderivative of the potential}, meaning by the potential
$q=s'\in W_{2,loc}^{-1}(\RR_+)$, where the derivative is understood in the sense of distribution theory.

By $\overline{l}$ we denote the differential expression constructed with $\overline{s}$ ---
the complex conjugate function to $s$. By $\overline{L}_0$ and $\overline{L}_M$ we denote the minimal and maximal operators,
defined by the differential expression $\overline{l}$. Their domains are, respectively, $\overline{\mathscr{D}}_0$
and $\overline{\mathscr{D}}$.

If the operator $L$ is given, $\widetilde{L}$ will denote its closure (if exists). The domain of $\widetilde{L}$ is denoted as usual by
$\mathscr{D}(\widetilde{L})$.

Everywhere $\|\cdot\|$ is the norm in $L_2(\RR_+)$. If $I$ is an interval (open or closed, or semiinterval), $\|\cdot\|_I$ is the norm in $L_2(I)$.

The length of the interval $I$ will be denoted by $|I|$.

By $C$ we denote arbitrary positive constants that may differ in adjacent formulas.

\begin{Statement}
\label{stker0}
The domain $\mathscr{D}_0$ is dense in $L_2(\RR_+)$, the operator $L_0$ is closable, its adjoint $L_0^*=\overline{L}_M$.
The kernel $\ker \widetilde{L_0}=\{0\}$, the defect $\widetilde{L_0}$ is finite: $\codim\im L_0\le 2$.
\end{Statement}

Due to finite defect, the resolvents of all extensions of $L_0$ (of all extensions that have resolvents) are either compact or not compact at the same time. Therefore,
all results below are formulated in terms of the minimal operator.

The following statement --- a necessary condition for the compactness of the resolvent --- is one of the key results of our study.
\begin{Theorem}
\label{thness}
If the minimal operator $L_0$ has extensions with compact resolvents, then
\begin{equation}
\label{eq2intthcrit}
\forall a>0\quad\lim\limits_{x\to+\infty}\int\limits_x^{x+a}\int\limits_x^{x+a}|s(\xi)-s(\eta)|^2\,d\xi d\eta=+\infty.
\end{equation}
\end{Theorem}

Under additional requirements on $s$, the condition \eqref{eq2intthcrit}
is a criterion for the compactness of the resolvent. Before these results, we give some definitions.

We will be interested in sectorial operators and quadratic forms (e.g., \cite[Ch. V,\S3]{Kato}).
\bigskip

\begin{Def} 
We call the operator $L$ (or the quadratic form $\mathfrak{l}$) defined on the domain $\mathscr{D}(L)$ (respectively, $\mathscr{D}(\mathfrak{l})$)
{\it sectorial} if the numerical range of the operator $L$ (of the form $\mathfrak{l}$) lies in some sector
$$
\Lambda_{\alpha,\beta}(z_0)=\set{z_0+z\in\CC}{\arg z\in[\alpha,\beta]},
$$
with $0\le\beta-\alpha\le\pi$, $z_0\in\CC$:
\begin{align*}
\forall& (y\in\mathscr{D}(L),\,\|y\|=1)\Rightarrow (Ly,y)\in\Lambda_{\alpha,\beta}(z_0), \mbox{ respectively,}\\
\forall& (y\in\mathscr{D}(\mathfrak{l}),\,\|y\|=1)\Rightarrow \mathfrak{l}[y]\in\Lambda_{\alpha,\beta}(z_0).
\end{align*}

In the case of strict inequality $0\le\beta-\alpha<\pi$, the operator (the form) is called {\it strictly sectorial}.

In the case $\alpha=\beta=0$ and $z_0\in\RR$, the operator (the form) is called {\it semi-bounded}.
\end{Def}
\bigskip

\begin{Def}
We say that $\sigma\in L_{p,unif}(\RR_+)$ for $p\ge1$, if $\sigma\in L_{p,loc}(\RR_+)$ and
$$
\sup_{\phantom{_+}x\in\RR_+}\int\limits_x^{x+1}|\sigma(\xi)|^p\,d\xi<+\infty.
$$
\end{Def}

\begin{Theorem}
\label{thcriter}
For $L_0$ to be strictly sectorial it is sufficient that there exist two complex-valued functions $f,g\in L_{2,unif}(\RR_+)$ such that for some $d>0$ and all $0\le x-t\le d$
\begin{equation}
\label{eqismagth2}
\arg \Bigl(s(x)+f(x)-s(t)-g(t)\Bigr)\in [\alpha,\beta],
\end{equation}
where $0\le\beta-\alpha<\pi$, $|\alpha|<\pi$, $|\beta|<\pi$.

In this case, the operator $L_0$ has extensions with compact resolvents if and only if the condition \eqref{eq2intthcrit} is satisfied.
\end{Theorem}

The following theorem is a direct generalization of Ismagilov's theorem \cite{Ismagilov}.
\begin{Theorem}
\label{thismagl}
Let $s$ satisfy the condition of the theorem \ref{thcriter}.
In this case, the operator $L_0$ has extensions with compact resolvents if and only if
\begin{equation}
\label{eqismag}
\forall A>0\ \forall \{\Delta_n\}_{n\in\NN}\ \lim_{n\to+\infty}\mu\set{(x,t)\in\Delta_n}{\bigl|s(x)-s(t)\bigr|<A}=0,
\end{equation}
where $\Delta_n=[a_n,a_n+d]\times [a_n,a_n+d]$, $a_n\to+\infty$ as $n\to\infty$ --- is a system of squares that goes to infinity; $\mu$ is the Lebesgue measure on $\RR^2$.
\end{Theorem}

The following theorem is a direct generalization of Brink's result \cite{Brinck}.
In the special case when $f,g=\const$ one can replace the limit in measure in \eqref{eqismag} by the usual one:
\begin{Theorem}
\label{thbrink}
For $L_0$ to be strictly sectorial it is sufficient that for
some $d>0$ and all $0\le x-t\le d$
$$
s(x)-s(t)\in\Lambda_{\alpha,\beta}(z_0),
$$
where $0\le\beta-\alpha<\pi$, $|\alpha|<\pi$, $|\beta|<\pi$, $z_0\in\CC$.

In this case, the operator $L_0$ has extensions with compact resolvents if and only if
$$
\forall a>0\ \lim\limits_{x\to+\infty}\bigl|s(x+a)-s(x)\bigr|=+\infty.
$$
\end{Theorem}

Natural questions arise:
\begin{itemize}
\item Is the condition \eqref{eqismagth2} of the theorem \ref{thcriter} complete in the description of the potentials of sectorial operators?

\item Are conditions of the type \eqref{eq2intthcrit} and \eqref{eqismag} sufficient for the compactness of the resolvent of a sectorial operator when refusing
restrictions of the type \eqref{eqismagth2}?
\end{itemize}

Zelenko \cite{Zelenko} gave an example of a semi-bounded operator satisfying the condition of the type \eqref{eqismag} with non-discrete spectrum.
We would like to look at this problem from a different angle. For simplicity, we turn to the real case ($\alpha=\beta=0$, $z_0\in\RR$).

Conditions \eqref{eq2intthcrit}, \eqref{eqismagth2}, \eqref{eqismag} are homogeneous in the following sense:
if $s$ satisfies any of them, 
then for arbitrary constant $k>0$ the function $ks$ will also satisfy it. At the same time,
the semi-boundedness itself is not a homogeneous property of the operator with respect to the potential:
if for $0<k\le1$ the semi-boundedness of the operator constructed by $s$ implies the semi-boundedness of the operator constructed by $ks$, 
then for $k>1$ this is not always true. Moreover,
\begin{Theorem}
\label{ths1s2exist}
There exist functions $s_1,s_2\in L_{2,loc}(\RR_+)$ and the number $\varepsilon>0$ such that the minimal operator $L_0$
\begin{itemize}
\item corresponding to $s_1$ is semi-bounded and has extentions with compact resolvents, while the operator corresponding to $(1+\varepsilon)s_1$ is not semi-bounded;
\item corresponding to $s_2$ is semi-bounded and has extentions with compact resolvents, while the operator corresponding to $(1+\varepsilon)s_2$ is semi-bounded but does not have 
extentions with compact resolvents.
\end{itemize}
\end{Theorem}
So the answers to all the questions posed above are negative.
The corresponding examples are constructed using the following Lemma, which involves the so-called Miura map \cite{Miura}:
\begin{Lemma}
\label{lmmiura}
The function $s\in L_{2,loc}(\RR_+)$ is a generalized antiderivative of the potential for some semi-bounded operator $L_0$ of the form \eqref{eqL0minop} if and only if
$$
s(x)=bx+\gamma(x)+\int\gamma^2(x)\,dx,
$$
for some real-valued $\gamma\in L_{2,loc}(\RR_+)$ and constant $b\in\RR$.
\end{Lemma}

Let us introduce one more function space. It is noteworthy that its definition involves the same integral as in theorems \ref{thness} and \ref{thcriter}.
\bigskip

\begin{Def}
We say that the function $q\in W_{2,unif}^{-1}(\RR_+)$ if $q\in W_{2,loc}^{-1}(\RR_+)$, and its antiderivative $s$ has the property:
\begin{equation}
\label{eqw2unifequivMain}
\sup_{\phantom{_+}x\in\RR_+}\int\limits_x^{x+1}\int\limits_x^{x+1}|s(\xi)-s(\eta)|^2\,d\xi d\eta<+\infty.
\end{equation}
\end{Def}

As follows from the lemma \ref{lmw2unif} of the Appendix, this definition of $W_{2,unif}^{-1}(\RR_+)$, is equivalent to the one in \cite{HrynivMykytyukPeriod}.

If the potential at least partially satisfies the homogeneity condition, its
perturbation by $\sigma'\in W_{2,unif}^{-1}(\RR_+)$ preserves both the sectorial property of the minimal operator
and its property to have extensions with compact resolvents:
\begin{Theorem}
\label{thunifpert}
Let $s\in L_{2,loc}(\RR_+)$, and for some $\varepsilon>0$
the minimal operator $L_0$ constructed by the antiderivative $(1+\varepsilon)s$ is strictly sectorial,
and the numerical range of $L_0$ lies in the sector $\Lambda_{\alpha,\beta}(z_0)$, where
$0\le\beta-\alpha<\pi$, $|\alpha|<\pi$, $|\beta|<\pi$, $z_0\in\CC$.

Then for any $\sigma\in L_{2,loc}(\RR_+)$ satisfying the condition \eqref{eqw2unifequivMain},
the minimal operator $L_{0,\sigma}$ constructed by $s+\sigma$ is also strictly sectorial.

If $L_0$ has extensions with compact resolvents,
then $L_{0,\sigma}$ also has extensions with compact resolvents.
\end{Theorem}

In the proofs of our theorems, the following lemma, is substantially used.

For an arbitrary interval $I=[a,b]$, we define the quadratic form $\mathfrak{l}_I$ in $L_2(I)$
on the domain $\mathscr{D}(\mathfrak{l}_I)=\overset{\circ}{W_2^1}(I)=\set{y\in W_2^1(I)}{y(a)=y(b)=0}$. For
$y\in \mathscr{D}(\mathfrak{l}_I)$:
$$
\mathfrak{l}_I[y]=\int\limits_I |y'(x)|^2\,dx - \int\limits_I s(x)\,d|y(x)|^2.
$$

\begin{Lemma}[Ismagilov's localization principle]
\label{lmlocalprinc}
For the operator $L_0$ to be sectorial it is necessary and sufficient that for some $d>0$ and any sequence of segments $I_k$ of fixed length
$|I_k|=d=\const>0$
there exists a common sector $\Lambda_{\alpha,\beta}(z_0)$, $0\le\beta-\alpha\le\pi$, $z_0\in\CC$ containing the numerical ranges of each of the forms
$\mathfrak{l}_k=\mathfrak{l}\vint{I_k}{}$:
$$
\forall k\in\NN\ \forall (y\in\mathscr{D}(\mathfrak{l}_k)\ \|y\|_{I_k}=1)\Rightarrow \mathfrak{l}_k[y]\in\Lambda_{\alpha,\beta}(z_0).
$$

In this case, the numerical range of the operator $L_0$ lies in some sector $\Lambda_{\alpha,\beta}(z_1)$ with the same angles $\alpha$ and $\beta$, but possibly with a different vertex
$z_1$.

For the compactness of the resolvents of extensions of the strictly sectorial ($\beta-\alpha<\pi$) operator $L_0$, it is necessary and sufficient that
\begin{equation}
\label{eqlimlemmaismag}
\lim\limits_{k\to\infty} \inf\limits_{
\begin{smallmatrix}y\in\mathscr{D}(\mathfrak{l}_k)\\
\|y\|_{I_k}=1
\end{smallmatrix}
}\bigl|\mathfrak{l}_k[y]\bigr|=+\infty,
\end{equation}
when the segments $I_k\to+\infty$, $|I_k|=d=\const>0$.
\end{Lemma}
\begin{Note}
For the sectoriality of $L_0$, as well as for the compactness of resolvents of extensions of the strictly sectorial $L_0$, 
it is sufficient to satisfy the corresponding conditions of the lemma \ref{lmlocalprinc} only for all segments of the form $[k-1,k]$ and $[k-1/2,k+1/2]$, $k\in\NN$.
\end{Note}
\bigskip

The further structure of the paper: section \ref{sect2} is entirely devoted to the proof of 
the necessary condition for the compactness of the resolvent (theorem \ref{thness}), 
section \ref{sect3} is devoted to the proof of Ismagilov’s localization principle (lemma \ref{lmlocalprinc}). In section \ref{sect4} the criteria for the 
compactness of the resolvent will be proved --- theorems \ref{thcriter}, \ref{thismagl}, \ref{thbrink} and \ref{thunifpert}. 
Section \ref{sect5} is devoted to the proof of the general form of the potential for the semi-bounded operator
and to the construction of counterexamples (theorem \ref{ths1s2exist}). 
And finally, section \ref{scappendix} (Appendix) contains auxiliary statements of minor interest.

\section{Necessary condition for the compactness of the resolvent}
\label{sect2}

{\noindent\bf Proof of theorem \ref{thness}.}
Assume the opposite, taking into account the inequality \eqref{equnifw2eq} of Appendix
we find the system $\{I_k\}_{k=1}^\infty$ of pairwise disjoint segments of equal length $|I_k|=\Delta>0$, the number $C>0$
and constants $\{C_k\}_{k=1}^\infty\subset\CC$ such that for any $k\in\NN$
$$
\int\limits_{I_k}|s(x)-C_k|^2\,dx<C.
$$
Without loss of generality, $\Delta<1/8$ and $C<1/8$. The latter can be easily achieved by reducing $I_k$.

For each $k\in\NN$, let $I_k=[a_k,b_k]$. Denote $s_k(x)=s(x+a_k)-C_k$, defining all $s_k$ on the common interval $x\in[0,\Delta]$.

Consider homogeneous systems of equations for $x\in[0,\Delta]$:
\begin{equation}
\label{eqsystemhomogen}
\begin{pmatrix}
y(x)\\
y^{[1]}(x)
\end{pmatrix}'=
\begin{pmatrix}
s_k(x) & 1 \\
-s_k^2(x) & -s_k(x)
\end{pmatrix}
\begin{pmatrix}
y(x)\\
y^{[1]}(x)
\end{pmatrix}=0
\end{equation}
and their fundamental solutions $(\varphi_k,\varphi_k^{[1]})$, $(\psi_k,\psi_k^{[1]})$, uniquely determined as the solution to the corresponding Cauchy problems
$(\varphi_k(0),\varphi_k^{[1]}(0))=(0,1)$, $(\psi_k(0),\psi_k^{[1]}(0))=(1,0)$.

Let us take any of these solutions, for example $(\varphi_k,\varphi_k^{[1]})$.

Let $m_k=\max|\varphi_k|$, $m_k^{[1]}=\max|\varphi_k^{[1]}|$, where the maxima are taken over the entire interval $[0,\Delta]$. 
Let $M_k=\max\{m_k,m_k^{[1]}\}$. 
It follows from \eqref{eqsystemhomogen}:
\begin{align*}
&m_k \le \int\limits_0^{\Delta}|s_k(x)|M_k\,dx+\int\limits_0^{\Delta}M_k\,dx \le M_k\Delta^{1/2}C^{1/2}+M_k\Delta,\\
&m_k^{[1]} \le 1+ \int\limits_0^{\Delta}|s_k(x)|^2M_k\,dx+\int\limits_0^{\Delta}|s_k(x)|M_k\,dx \le1+ M_kC+M_k\Delta^{1/2}C^{1/2},
\end{align*}
obviously
$$
M_k\le 1+2M_k\Delta^{1/2}C^{1/2}+M_k\Delta+M_kC=1+M_k(\Delta^{1/2}+C^{1/2})^2\le1+M_k/2.
$$
The last inequality follows from the condition $\Delta<1/8$ and $C<1/8$. Finally, $M_k\le 2$, and this estimate is uniform in $k\in\NN$, 
and is valid for both solutions $(\varphi_k,\varphi_k^{[1]})$
and $(\psi_k,\psi_k^{[1]})$.

Further, we show that each of the systems $\{\varphi_k\}_{k=1}^\infty$ and $\{\psi_k\}_{k=1}^\infty$ satisfies the conditions of the Arzela-Ascoli theorem.

We have just shown the uniform boundedness, the equicontinuity follows from the inequality below with $0\le x_1<x_2\le\Delta$:
\begin{gather*}
|\varphi_k(x_2)-\varphi_k(x_1)|\le\int\limits_{x_1}^{x_2}|s_k(x)\varphi_k(x)+\varphi_k^{[1]}(x)|\,dx\le\\
\le2(x_2-x_1)^{1/2}\Bigl(\int\limits_{x_1}^{x_2}|s_k(x)|^2\,dx\Bigr)^{1/2}+2(x_2-x_1)\
\le2(x_2-x_1)^{1/2}C^{1/2}+2(x_2-x_1),
\end{gather*}
where we used the estimate $M_k\le 2$. A similar inequality is also valid for $\psi_k$.

The Arzela-Ascoli theorem implies the possibility of extracting uniformly convergent 
subsequences from $\{\varphi_k\}_{k=1}^\infty$ and $\{\psi_k\}_{k=1}^\infty$. Without loss of generality, we assume that each of initial 
sequences is uniformly convergent (otherwise we restrict ourselves to considering only those $I_k$ for which uniform convergence takes place):
\begin{equation}
\label{eqvarphipsiunif}
\varphi_k(x)\rightrightarrows\varphi(x),\quad \psi_k(x)\rightrightarrows\psi(x),\quad x\in[0,\Delta].
\end{equation}

Now we show the linear independence of $\varphi$ and $\psi$ on any segment $[0,a]\subset[0,\Delta]$. Assume, on the contrary, that
there exist two constants $|C_1|+|C_2|>0$ and
$C_1\varphi(x)+C_2\psi(x)\equiv 0$ for all $x\in[0,a]$.

We denote $\chi_k(x)=C_1\varphi_k(x)+C_2\psi_k(x)$, $\chi_k(x)\rightrightarrows0$ for $x\in[0,a]$. 

We introduce
$$
A_k(x)=\int\limits_0^xs_k(\xi)\chi_k^{[1]}(\xi)\,d\xi,
$$
the functions $\{A_k\}_{k=1}^\infty$ satisfy the conditions of the Arzela-Ascoli theorem, we assume, without loss of generality,
the existence of the uniform limit:
$$
A_k(x)\rightrightarrows A(x),\quad x\in[0,\Delta].
$$

Integrating \eqref{eqsystemhomogen} over the arbitrary interval $[0,x]\subset[0,a]$, substituting $\chi_k$ as a solutions, we obtain:
$$
\chi_k(x)=\chi_k(0)+\chi_k^{[1]}(0)x+\int\limits_0^xs_k(\eta)\chi_k(\eta)\,d\eta-\int\limits_0^x\int\limits_0^\xi s_k^2(\eta)\chi_k(\eta)\,d\eta\, d\xi-
\int\limits_0^x A_k(\eta)\,d\eta.
$$

Since $\chi_k(0)=C_2\rightrightarrows0$, we conclude that $C_2=0$. Further, we calculate $\chi_k^{[1]}(0)=C_1$, 
and take the uniform limit as $k\to\infty$ in the last equality, and finally obtain:
$$
0=C_1 x-\int\limits_0^x A(\eta)\,d\eta,
$$
whence $A(x)\equiv C_1$. Since for all $k\in\NN$ $A_k(0)=0$, then $A(0)=C_1=0$, we come to a contradiction with the linear dependence of
$\varphi$ and $\psi$ for $x\in[0,a]$.

In order to come to a contradiction with the main assumption of the theorem, we will prove the existence of the sequence of functions $\{y_k\}_{k=1}^\infty\subset \mathscr{D}_{L_0}$,
such that
for some $S>0$ for all $k\in\NN$
$\supp y_k\subset I_k$, $\|y_k\|=1$, $\|L_0y_k\|<S$.

Instead of considering $L_0$ in $L_2(I_k)$ for each $k\in\NN$, we introduce operators $L_k$ in $L_2[0,\Delta]$ obtained
from $L_0$ by shifting of the independent variable by $a_k$.
It is essential that adding the constant to $s$ does not change the operator itself. The corresponding differential expression obtained by the shift is denoted by $l_k$.

It is sufficient to prove the existence of $\{y_k\}_{k=1}^\infty\subset \mathscr{D}_{L_k}$ for some $S>0$,
$\supp y_k\subset [0,\Delta]$, $\|y_k\|=1$, $\|L_ky_k\|<S$.

For arbitrary $k\in\NN$ and $f\in L_2[0,\Delta]$ we set
$$
(R_kf)(x)=\psi_k(x)\int\limits_0^x\varphi_k(\xi)f(\xi)\,\xi-\varphi_k(x)\int\limits_0^x\psi_k(\xi)f(\xi)\,\xi,
$$
as is easy to verify, $R_k$ is a Cauchy operator and $y=R_kf$ is a solution to the equation $l_ky=f$ with initial conditions $y(0)=y^{[1]}(0)=0$.

Let $\overset{\circ}R_k$ be the restriction of $R_k$ to the space $f\in L_2[0,\Delta]$, $f\perp\langle \overline\varphi_k,\overline\psi_k\rangle$, where
$\langle\ldots\rangle$ denotes the linear span. The existence of $\{y_k\}_{k=1}^\infty$ satisfying the above mentioned properties is equivalent to 
the uniform separability of the norm $\overset{\circ}R_k$ from 0.

By the uniform convergence in \eqref{eqvarphipsiunif}, uniformly $R_k\rightrightarrows R$ to some compact operator $R$, which, by lemma
\cite[L.2]{TumMatSbor24}, is not finite-dimensional, at least its singular value $s_3>0$.

In this case there exists $a>0$ such that for all $k\in\NN$ the corresponding singular values of all operators $s_3(R_k)>a$. Let us evaluate the norms of $\overset{\circ}R_k$:
$$
\|\overset{\circ}{R}_k\|^2=\sup_{f\perp\langle\overline\varphi_k,\overline\psi_k\rangle}
\frac{\|R_kf\|^2}{\|f\|^2}\ge \min_{\mathscr{L},\,\dim\mathscr{L}=2}\ \max_{f\perp\mathscr{L}}
\frac{\bigl((R_k)^*R_kf,f\bigr)}{\|f\|^2}=s_3^2(R_k)>a^2\,
$$
the minimum is taken over all two-dimensional subspaces in $L_2[0,\Delta]$. The equality follows from \cite[Ch.II,\S1]{GhbKrn}, which completes the proof. \qquad$\Box$
\section{Localization principle}
\label{sect3}
We start with some auxiliary statements and finally prove the lemma \ref{lmlocalprinc}. 

\begin{Lemma}
\label{lmnormcrit}
The minimal operator $L_0$ of the form \eqref{eqL0minop} has extensions
with compact resolvents if and only if the set
$$
\mathscr{Y}=\set{y\in \mathscr{D}_0}{\|L_0y\|+\|y\|\le1}
$$
is precompact in $L_2(\RR_+)$.
\end{Lemma}

{\noindent\bf Proof.} Suppose $L_0$ has the extension $L$ with compact resolvent. We take $0\ne\lambda_0\in\CC$,
such that $R_{\lambda_0}=(L-\lambda_0)^{-1}$ is a compact operator. For any $y\in\mathscr{Y}$
$$
\|(L-\lambda_0)y\|=\|(L_0-\lambda_0)y\|\le\|L_0y\|+|\lambda_0|\|y\|\le C(\|L_0y\|+\|y\|)\le C
$$
--- is bounded, hence $\mathscr{Y}$ is precompact.

Now we prove the converse statement. Let $\mathscr{Y}$ be precompact. Consider $\widetilde{L_0}$ --- the closure of $L_0$. The following set is also precompact:
$$
\mathscr{X}=\set{x\in \mathscr{D}(\widetilde{L_0})}{\|\widetilde{L_0}x\|+\|x\|\le1}.
$$

Let us show that the image of $\widetilde{L_0}$ is closed. Let $f_n=\widetilde{L_0}x_n\to f$ as $n\to\infty$. The sequence $x_n$ contains a bounded
subsequence. If this were not so, then $\|x_n\|\to+\infty$, and $\widetilde{L_0}x_n/\|x_n\|\to0$, and for large $n>n_0$,
$x_n/\|x_n\|/2\in\mathscr{X}$, extracting the subsequence $x_{n_l}/\|x_{n_l}\|\to\xi$, $\|\xi\|=1$, $\widetilde{L_0}\xi=0$, 
i.e. $\ker\widetilde{L_0}\ne\{0\}$, which contradicts the statement \ref{stker0}.

Thus, some subsequence $x_{n_l}$ is bounded. We can extract the convergent subsequence $x_{n_{l_m}}\to x_0$ due to the precompactness 
of $\mathscr{X}$. Since $\widetilde{L_0}$ is closed, $\widetilde{L_0}x_0=f$, which proves that the image of $\widetilde{L_0}$ is closed.

Thus $\mathfrak{M}=\im \widetilde{L_0}$ is a closed space, $\widetilde{L_0}$ maps $\mathscr{D}(\widetilde{L_0})$ one-to-one onto
$\mathfrak{M}$, and the inverse operator $R:\mathfrak{M}\to\mathscr{D}(\widetilde{L_0})$ is compact due to the precompactness of $\mathscr{X}$.

If $\mathfrak{M}=L_2(\RR_+)$, then $\widetilde{L_0}$ is the required extension of $L_0$ with the compact resolvent.

Otherwise we take the basis $f_1,\ldots,f_m$ ($1\le m\le2$) in $\mathfrak{N}=L_2(\RR_+)\ominus\mathfrak{M}$ and the same number of linearly independent vectors
over $\mathscr{D}(\widetilde{L_0})$: $y_1,\ldots,y_m$ and extend $R$ to $\mathfrak{N}$, setting $Rf_j=y_j$ on the basis. The resulting operator is compact and
injective; its inverse $L$ on the domain $\mathscr{D}=\mathscr{D}(\widetilde{L_0})+\langle y_1,\ldots,y_m\rangle$ --- is the required extension of $L_0$.\qquad$\Box$
\bigskip

Studying the compactness of resolvents for positive operators the following lemma plays a crucial role:
\begin{LemmaNoCounter}[Rellich]
Let $A$ be the self-adjoint operator in the Hilbert space $\mathfrak{H}$ defined on the domain $\mathscr{D}(A)$ satisfying the condition
$$
(Ay,y)\ge(y,y)
$$
for all $y\in\mathscr{D}(A)$. The resolvent of $A$ is compact if and only if the set of all $y\in\mathscr{D}(A)$ such that
$$
(Ay,y)\le1,
$$
is precompact in $\mathfrak{H}$.
\end{LemmaNoCounter}

We will need its generalization:

\begin{Lemma}
\label{lmrlchtum}
The minimal strictly sectorial operator $L_0$ of the form \eqref{eqL0minop} has extensions
with compact resolvents if and only if the set
$$
\mathscr{Y}=\set{y\in \mathscr{D}_0}{|(L_0y,y)|+\|y\|^2\le1}
$$
is precompact in $L_2(\RR_+)$.
\end{Lemma}

{\noindent\bf Proof.} We find $\gamma\in\RR$ and $z_0\in\CC$ such that the numerical range of $M_0=e^{i\gamma}(L_0-z_0I)$ lies
in the sector $|\arg(z-z_0)|\le\theta<\pi/2$.
The sesquilinear form defined by $M_0$ is closable \cite[Ch.VI,\S1,Th.1.27]{Kato}. We denote its closure by $\mathfrak{m}$. Let $\mathscr{D}(\mathfrak{m})$
be its domain. The domain $\mathscr{D}_0$ is a core of $\mathfrak{m}$ in the sense of \cite[Ch.VI,\S1]{Kato}.

By the first representation theorem
\cite[Ch.VI,\S2,Th.2.1]{Kato}, $\mathfrak{m}$ is associated with the closed operator $M_F\supset M_0$ --- the Friedrichs extension of $M_0$
\cite[Ch.VI,\S2,p.3]{Kato}.
In view of the finite codimension
of the image of $M_0$, the existence of extensions of $L_0$ with compact resolvents is equivalent to the existence of a compact resolvent of $M_F$.

We construct the symmetric form $\mathfrak{t}=\re\mathfrak{m}$: for $u,v\in\mathscr{D}(\mathfrak{m})$,
we set
$$
\mathfrak{t}[u,v]=\frac{1}{2}(\mathfrak{m}[u,v]+\overline{\mathfrak{m}[v,u]}).
$$

The form $\mathfrak{t}$ is closed on the domain $\mathscr{D}(\mathfrak{t})=\mathscr{D}(\mathfrak{m})$, densely defined, and non-negative;
by the second representation theorem \cite[Ch.VI,\S2,Th.2.23]{Kato} it is associated with the self-adjoint operator $T\ge0$,
$\mathscr{D}(\mathfrak{t})=\mathscr{D}(T^{1/2})\supset\mathscr{D}(T)$, where the domain $\mathscr{D}(T)$ of $T$ is the core of $\mathfrak{t}$, and
$$
\mathfrak{t}[u,v]=(T^{1/2}u,T^{1/2}v),\quad u,v\in\mathscr{D}(\mathfrak{t}).
$$
The domain $\mathscr{D}_0$ is also the core of $\mathfrak{t}$.

The resolvents of both operators $M_F$ and $T$ are compact or not compact at the same time \cite[Ch.VI,\S3,Th.3.3]{Kato}.

The following arguments use the Rellich's Lemma for $T+I$.

Let us prove the necessity. Let the resolvent of $M_F$ be compact. For any $y\in\mathscr{Y}$,
$$
\mathfrak{t}[y,y]=\re(M_0y,y)\le|(M_0y,y)|\le|(L_0y,y)|+|z_0|\|y\|^2\le 1+|z_0|=C_0.
$$
Since $\mathscr{D}(T)$ is a core of $\mathfrak{t}$,
for any $\varepsilon>0$ and any $y\in\mathscr{Y}$ we find $u\in\mathscr{D}(T)$ such that
$|\mathfrak{t}[u,u]-\mathfrak{t}[y,y]|<\varepsilon$ and $\|u-y\|<\varepsilon$. The compactness of the resolvent of
$M_F$ implies the compactness of the resolvent of $T+I$. Since
$$
((T+I)u,u)=\mathfrak{t}[u,u]+\|u\|^2\le (\mathfrak{t}[y,y]+\varepsilon) + (\|y\|+\varepsilon)^2
\le (C_0+\varepsilon) + (2+2\varepsilon^2),
$$
from Rellich's lemma we conclude that the entire set $\{u\}$ is precompact, and in view of the arbitrariness of $\varepsilon$
-- we conclude that $\mathscr{Y}$ is precompact.

Let us prove sufficiency.
Let $\mathscr{Y}$ be precompact. 

Consider $U=\set{u\in\mathscr{D}(T)}{((T+I)u,u)\le1}$.
Since $\mathscr{D}_0$ is a core of $\mathfrak{t}$,
for any $\varepsilon>0$ and any $u\in U$ we find $y\in\mathscr{D}_0$ such that
$|\mathfrak{t}[u,u]-\mathfrak{t}[y,y]|<\varepsilon$ and $\|u-y\|<\varepsilon$. We denote the set of all such $y$
by $\mathscr{Y}_\varepsilon$.

Due to the strict sectoriality of $M_0$,
the following estimate is valid for any $y\in\mathscr{Y}_\varepsilon$:
$$
|(M_0y,y)|\le C_1\re(M_0y,y)=C_1\mathfrak{t}[y,y]\le C_1(\mathfrak{t}[u,u]+\varepsilon)
\le C_1(1+\varepsilon)
$$
for some $C_1>0$ depending
only on $\theta<\pi/2$ --- the half-angle of the sector. Since
\begin{gather*}
|(L_0y,y)|+\|y\|^2\le|(M_0y,y)|+(|z_0|+1)\|y\|^2\le\\
\le C_1\re(M_0y,y)+(|z_0|+1)(\varepsilon+\|u\|)^2\le
C_1(\mathfrak{t}[u,u]+\varepsilon)+(|z_0|+1)(2+2\varepsilon^2)\le\\
\le C_1(1+\varepsilon)+(|z_0|+1)(2+2\varepsilon^2)=C_2,
\end{gather*}
the set $\mathscr{Y}_\varepsilon$ is precompact. Since $\varepsilon>0$ is arbitrary, $U$ is precompact.
Applying Rellich's lemma, we conclude that $T$ has a compact resolvent,
therefore so does $M_F$.\qquad$\Box$
\bigskip

It is convenient to prove the lemma \ref{lmlocalprinc} in the language of quadratic forms.

If $y,y^{[1]}\in AC_{loc}(\RR_+)$ and $s\in L_{2,loc}(\RR_+)$, then $y'\in L_{2,loc}(\RR_+)$. Consider 
for $u,v\in \mathscr{D}_0$:
\begin{equation}
\label{eqformmS}
\mathfrak{l}_0[u]=(L_0 u,u)=\int\limits_0^{+\infty} |u'|^2\,dx - \int\limits_0^{+\infty} s\,d|u|^2.
\end{equation}

The right-hand side of \eqref{eqformmS} allows us to extend the form $\mathfrak{l}_0\subset\mathfrak{l}$ to the domain
\begin{equation}
\label{eqformminD}
\mathscr{D}(\mathfrak{l})=\set{y\in W_{2,loc}^1(\RR_+)}{y(0)=0,\, \supp y\mbox{ --- is a compact}},
\end{equation}

\begin{Def}
We call the quadratic form $\mathfrak{l}$ defined by the right-hand side of \eqref{eqformmS} on the domain \eqref{eqformminD} the {\it minimal
quadratic form} defined by the expression \eqref{eqformmS}.
\end{Def}
\bigskip

Despite that $\mathscr{D}_0\subset\mathscr{D}(\mathfrak{l})$ in the strict sense, this definition is natural because of the following statement:
\begin{Statement}
\label{stsect}
If the form $\mathfrak{l}_0[u]=(L_0u,u)$ is closable, then $\widetilde{\mathfrak{l}}_0=\widetilde{\mathfrak{l}}$.

The operator $L_0$ is sectorial if and only if the form $\mathfrak{l}$ is sectorial. The sectors containing the numerical ranges of the
sectorial operator $L_0$ and the form $\mathfrak{l}$ coincide.
\end{Statement}
{\noindent\bf Proof.} Due to lemma \ref{lmDdenseW21} of Appendix for any $u\in\mathscr{D}(\mathfrak{l})$ with $\supp u\subset [0,a]$
exists the sequence $u_n\in\mathscr{D}_0$,
$u_n\to u$ in the metric of $W_{2}^1[0,a]$, thus $\mathfrak{l}_0[u_n]\to \mathfrak{l}[u]$.\qquad$\Box$
\bigskip

Below we use identical notations for the quadratic form $\mathfrak{l}[\cdot]$ and the sesquilinear form $\mathfrak{l}[\cdot,\cdot]$,
where for $u,v\in\mathscr{D}(\mathfrak{l})$ \eqref{eqformminD}:
$$
\mathfrak{l}[u,v]=\int\limits_0^{+\infty} u'\,\overline{v}'\,dx - \int\limits_0^{+\infty} s\,d(u\,\overline{v}).
$$
\begin{Lemma}
\label{lmDeltasum}
Let $\Delta_k\subset\RR_+$, $k\in\NN$, be intervals (open, closed or semi-open, possibly unbounded) with the properties:
\begin{itemize}
\item union of all $\cup\Delta_k=\RR_+$;
\item pairwise intersection of adjacent $\Delta_k\cap\Delta_{k+1}=\delta_k\ne\varnothing$, $k\in\NN$;
\item all other intersections $\Delta_k\cap\Delta_{k+j}=\varnothing$ for $j>1$, $k\in\NN$.
\end{itemize}

Let for all $k\in\NN$ the real-valued functions $\varphi_k\in W_{2,loc}^1(\RR_+)$, $\supp\varphi_k\subset \overline\Delta_k$ (in closure),
and their squares form a partition of unity:
$$
\sum_{k=1}^{\infty}\varphi_k^2(x)\equiv 1.
$$
For any $y\in\mathscr{D}_\mathfrak{l}$ we set $y_k=\varphi_k^2y$, $u_k=\varphi_k\varphi_{k+1}y$.

The following identity is valid:
\begin{equation}
\label{eqismagilovsum}
\mathfrak{l}[y]=\sum_{k=1}^{\infty}\mathfrak{l}[y_k]+2\sum_{k=1}^{\infty}\mathfrak{l}[u_k]-
2\sum_{k=1}^{\infty}\int\limits_0^{+\infty}(\varphi_k'\varphi_{k+1}-\varphi_k\varphi_{k+1}')^2|y|^2\,dx.
\end{equation}
\end{Lemma}
{\noindent\bf Proof.} Since $\supp y$ is compact, for some $N\in\NN$, all $y_k\equiv0$ as $k>N$. Then
$$
\mathfrak{l}[y]=\sum_{k,l=1}^{N}\mathfrak{l}[y_k,y_l]=\sum_{k=1}^{N}\mathfrak{l}[y_k]+\sum_{k=1}^{N-1}\bigl(\mathfrak{l}[y_k,y_{k+1}]+\mathfrak{l}[y_{k+1},y_k]\bigr).
$$

For any $k\in\NN$:
\begin{gather*}
2\mathfrak{l}[u_k]-\Bigl(\mathfrak{l}[y_k,y_{k+1}]+\mathfrak{l}[y_{k+1},y_k ]\Bigr)=2\int\limits_0^{+\infty}\bigl(|u_k'|^2-\re(y_k'\overline{y_{k+1}'})\bigr)\,dx=\\
=2\int\limits_ 0^{+\infty}\bigl([(\varphi_k\varphi_{k+1})']^2|y|^2-(\varphi_k^2)'(\varphi_{k+1}^2)'|y|^2\bigr)\,dx=
2\int\limits_0^{+\infty}(\varphi_k'\varphi_{k+1}-\varphi_k\varphi_{k+1}')^2|y|^2\,dx,
\end{gather*}
which completes the proof.\qquad$\Box$
\bigskip

{\noindent\bf Proof of the localization principle (lemma \ref{lmlocalprinc}).} Let us start with sectoriality. Only sufficiency is of interest.
Let $U_k=[k-1,k]$ and $V_k=[k-1/2,k+1/2]$, $k\in\NN$. The aggregate system $\sint{U_l}{l\in\NN}{}\cup \sint{V_m}{m\in\NN}{}$, 
reordered in the natural way, forms the system $\sint{\Delta_k}{k\in\NN}{}$ satisfying the conditions of Lemma \ref{lmDeltasum}.

Denote for any $k\in\NN$, $\delta_k=\Delta_k\cap\Delta_{k+1}$ --- intervals of length $1/2$.

For each $\Delta_k=[a_k,a_k+1]$, we set
$$
\varphi_k(x)=
\left\{
\begin{array}{cl}
\sin(\pi(x-a_k))&\mbox{ for } x\in \Delta_k\\
0&\mbox{ for } x\not\in\Delta_k
\end{array}
\right.,
$$
the functions $\varphi_k$ satisfy the conditions of lemma \ref{lmDeltasum}.

We take an arbitrary $y\in\mathscr{D}(\mathfrak{l})$, $\|y\|=1$, let $y_k=\varphi_k^2y$, $u_k=\varphi_k\varphi_{k+1}y$.

By the conditions of our lemma,
$$
\mbox{for }y_k\not\equiv0,\ \frac{\mathfrak{l}[y_k]}{\|y_k\|^2}\in \Lambda_{\alpha,\beta}(z_0);\quad
\mbox{for }u_k\not\equiv0,\ \frac{\mathfrak{l}[u_k]}{\|u_k\|^2}\in \Lambda_{\alpha,\beta}(z_0).
$$

Keeping in mind the identity \eqref{eqismagilovsum}, we turn to the following sum, taking into account the finite number of its terms:
\begin{gather*}
\sum_{k=1}^{\infty}\mathfrak{l}[y_k]+2\sum_{k=1}^{\infty}\mathfrak{l}[u_k]=\\
=\sum_{
\begin{smallmatrix}k=1\\
y_k\not\equiv0
\end{smallmatrix}
}^{\infty}\frac{\mathfrak{l }[y_k]}{\|y_k\|^2}\,\|y_k\|^2+\sum_{\begin{smallmatrix}k=1\\
u_k\not\equiv0
\end{smallmatrix}}^{\infty}\frac{\mathfrak{l}[u_k]}{\|u_k\|^2}\,2\|u_k\|^2\in\Lambda_{\alpha,\beta}(z_0),
\end{gather*}
since
$$
1=\|y\|^2=\sum_{
\begin{smallmatrix}k=1\\
y_k\not\equiv0
\end{smallmatrix}
}^{\infty}\|y_k\|^2+\sum_{\begin{smallmatrix}k=1\\
u_k\not\equiv0
\end{smallmatrix}}^{\infty}2\|u_k\|^2,
$$
and $\Lambda_{\alpha,\beta}(z_0)$ is a convex set.

Let's evaluate the last sum in \eqref{eqismagilovsum}:
\begin{gather}
\notag
2\sum_{k=1}^{\infty}\int\limits_0^{+\infty}(\varphi_k'\varphi_{k+1}-\varphi_k\varphi_{k+1}')^2|y|^2\,dx=\\
=2\sum_{k=1}^{\infty}\int\limits_{\delta_k}(\varphi_k'\varphi_{k+1}-\varphi_k\varphi_{k+1}')^2|y|^2\,dx=
2\pi^2\sum_{k=1}^{\infty}\int\limits_{\delta_k}|y|^2\le2\pi^2.
\label{eqlstsum2p2}
\end{gather}

The statement about sectoriality follows from the lemma \ref{lmDeltasum} and simple geometric reasoning.

Now we turn to the second part of the theorem and prove its sufficiency. Assume the contrary, that \eqref{eqlimlemmaismag} holds, but
$L_0$ has no extensions with compact resolvents.

Consider the set $\mathscr{Y}$ provided by the lemma \ref{lmnormcrit}. It is not precompact. For any
$T>0$ we construct a partition of unity on $[T,T+1]$, setting for $x\in[T,T+1]$:
$$
\tau_T(x)=30\int\limits_x^{T+1}(\xi-T)^2(\xi-T-1)^2\,d\xi,\quad\chi_T(x)=1-\tau_T(x).
$$
On $[T,T+1]$ the function $\tau_T$ is monotonically decreasing, $\tau_T(T)=1$, $\tau_T(T+1)=\tau_T'(T+1)=\tau_T'(T)=0$.
We extend both functions $\tau_T$ and $\chi_T$ by constants $0$ and $1$ by continuity to the whole $\RR_+$.

Each $y\in\mathscr{Y}$ can be represented as $y=y_{T,0}+y_{T,\infty}$, where $y_{T,0}=y\cdot\tau_T$,
$y_{T,\infty}=y\cdot\chi_T$. The set $\{y\vint{[0,T+1]}{}\}$ of cuts of $\mathscr{Y}$ on $[0,T+1]$ is precompact in
$L_2(0,T+1)$, which follows from the compactness of the Cauchy operator in $L_2(0,T+1)$. The latter
associates with any $f\in L_2(0,T+1)$ a solution to the Cauchy problem $l(y)=f$ with zero initial conditions $y(0)=y^{[1]}(0)=0$.
The Cauchy operator itself has an explicit form:
$$
(Rf)(x)=\psi(x)\int\limits_0^x\varphi(\xi)f(\xi)\,\xi-\varphi(x)\int\limits_0^x\psi(\xi)f(\xi)\,\xi,
$$
where $\varphi,\psi\in AC[0,T+1]$ are independent solutions to the equation $l(y)=0$ with initial conditions
$(\varphi(0),\varphi^{[1]}(0))=(0,1)$, $(\psi(0),\psi^{[1]}(0))=(1,0)$.

Since $|\tau_T|<1$, the set $\{y_{T,0}\}$ is precompact in $L_2(\RR_+)$.

There exists $\varepsilon_0>0$ such that for any $T>0$ we can find the element $y_T^*\in\mathscr{Y}$ with $\supp y_T^*\subset[T,+\infty)$, 
$\|y_T^*\|>\varepsilon_0$. At least one of the elements of $\{y_{T,\infty}\}$ has such property. Otherwise, for any
$\varepsilon>0$ we would find $T_0>0$ for which all $|y_{T_0,\infty}|<\varepsilon$. Given the precompactness of $\{y_{T_0,0}\}$, there exists a finite $\varepsilon$-net 
of elements $y_{T_0,0}^1,\ldots,y_{T_0,0}^k$. 
They, being themselves elements of $\mathscr{Y}$, form a $2\varepsilon$-net in it, which is impossible, since $\mathscr{Y}$ is not precompact.

Let's take the same system of intervals $\Delta_k$ as in the first part of the proof.
Further we set $\theta=(\beta-\alpha)/2$, $0\le\theta<\pi/2$.

We take sufficiently large $A$:
\begin{equation}
\label{eqAsec2p2}
A>\Bigl(
\frac{1}{\varepsilon_0^2}+|z_0|+2\pi^2
\Bigr)\sec\theta
\end{equation}
and find $N\in\NN$ such that for all $k\ge N$
\begin{equation}
\label{eqlkz0A}
\inf\limits_{
\begin{smallmatrix}y\in\overset{\circ}{W_2^1}(\Delta_k)\\
\|y\|_{\Delta_k}=1
\end{smallmatrix}
}\bigl|\mathfrak{l}_k[y]-z_0\bigr|>A.
\end{equation}

We take the element $y_N^*\in\mathscr{D}_0$ so that $\supp y_N^*\subset[N,+\infty)$, $\|y_N^*\|>\varepsilon_0$ and 
\begin{equation}
\label{eqyNasterisk}
\|L_0y_N^*\|+\|y_N^*\|\le1.
\end{equation}

As in the first part of the proof, we denote $y_{N,k}^*=\varphi_k^2y_N^*$, $u_{N,k}^*=\varphi_k\varphi_{k+1}y_N^*$.
By our assumption,
$$
\mbox{for }y_{N,k}^*\not\equiv0,\ \frac{\mathfrak{l}[y_{N,k}^*]}{\|y_{N,k}^*\|^2}\in \Lambda_{\alpha,\beta}(z_0);\quad
\mbox{for }u_{N,k}^*\not\equiv0,\ \frac{\mathfrak{l}[u_{N,k}^*]}{\|u_{N,k}^*\|^2}\in \Lambda_{\alpha,\beta}(z_0).
$$

We introduce a right orthogonal coordinate system $(\xi,\eta)$ on the plane, placing $z_0$ in its center, directing the $\xi$ axis along the bisector of $\Lambda_{\alpha,\beta}(z_0)$.
In the introduced coordinate system $(\xi,\eta)\in\Lambda_{\alpha,\beta}(z_0)$ if and only if $\xi\ge0$, $|\eta|\le \xi\tan\theta$.

It follows from the condition \eqref{eqlkz0A} that for all $k\ge N$:
$$
\mbox{for }y_{N,k}^*\not\equiv0,\ \frac{\mathfrak{l}[y_{N,k}^*]}{\|y_{N,k}^*\|^2}\in \Lambda_{\alpha,\beta}^A(z_0);\quad
\mbox{for }u_{N,k}^*\not\equiv0,\ \frac{\mathfrak{l}[u_{N,k}^*]}{\|u_{N,k}^*\|^2}\in \Lambda_{\alpha,\beta}^A(z_0),
$$
where $\Lambda_{\alpha,\beta}^A(z_0)\subset \Lambda_{\alpha,\beta}(z_0)$ is the set of $(\xi,\eta)\in\Lambda_{\alpha,\beta}(z_0)$ for which $\xi\ge A\cos\theta$. 
This set is a convex one.

Using the convexity of $\Lambda_{\alpha,\beta}^A(z_0)$, we obtain:
$$
\frac{1}{\|y_N^*\|^2}\Bigl(
\sum_{k=1}^{\infty}\mathfrak{l}[y_{N,k}^*]+2\sum_{k=1}^{\infty}\mathfrak{l}[u_{N,k}^*]
\Bigr)
\in \Lambda_{\alpha,\beta}^A(z_0),
$$
and, therefore, taking into account \eqref{eqismagilovsum} and \eqref{eqlstsum2p2},
\begin{equation}
\label{eqL0vAcosth}
\Bigl|(L_0\frac{y_N^*}{\|y_N^*\|},\frac{y_N^*}{\|y_N^*\|})-z_0\Bigr|\ge A\cos\theta-2\pi^2.
\end{equation}

But in view of \eqref{eqyNasterisk}
$$
|(L_0y_N^*,y_N^*)|\le\|L_0y_N^*\|\|y_N^*\|\le1,
$$
and since $\|y_N^*\|>\varepsilon_0$,
$$
\Bigl|(L_0\frac{y_N^*}{\|y_N^*\|},\frac{y_N^*}{\|y_N^*\|})-z_0\Bigr|<
\frac{1}{\varepsilon_0^2}+|z_0|,
$$
which contradicts \eqref{eqAsec2p2} and \eqref{eqL0vAcosth}.

Now we prove the necessity. Assume, on the contrary, that $L_0$ has an extension with compact resolvent, but
there exists a sequence of pairwise disjoint segments of equal length $I_k$ and $y_k\in\overset{\circ}{W_2^1}(I_k)$,
$\|y_k\|=1$, and $|\mathfrak{l}_k[y_k]|<C$. Taking into account lemma \ref{lmDdenseW21} of the Appendix, we find the sequence $u_k\in\mathscr{D}_0$,
$\supp u_k\subset I_k$, $\|u_k\|=1$, $|(L_0u_k,u_k)|<2C$.

The sequence $u_k$ is obviously not precompact, thus by lemma \ref{lmrlchtum} the operator $L_0$ cannot have extensions with compact resolvents,
we have arrived at a contradiction. The lemma is completely proved.\qquad$\Box$

\section{Criteria for compactness of the resolvent}
\label{sect4}

In this section we prove theorems \ref{thcriter}, \ref{thismagl}, \ref{thbrink} and \ref{thunifpert}.
\bigskip

{\noindent\bf Proof of theorem \ref{thcriter}.} Denote the closed sector $\Lambda=\{\arg z\in [\alpha,\beta]\}$.
Without loss of generality, either
$\alpha=\beta=0$ or $-\pi<\alpha<0<\beta<\pi$, (otherwise, $\Lambda$ should be expanded, while remaining within the conditions of the theorem).

We introduce a right-handed coordinate system $(u,v)$ on the plane, placing its center at $0$, and
directing the $u$ axis in the direction of the ray $\arg z=\alpha$, and the $v$ axis --- in the direction of $\arg z=\beta$.
The projection onto the corresponding axis parallel to the second one will be denoted by the subscript $u$ or $v$.

Thus, for any $z\in\CC$, $z=e^{i\alpha}z_u+e^{i\beta}z_v$. Simplifying
the notation, we write $z=(z_u,z_v)$. Respectively
$z=(z_u,z_v)\in\Lambda$ if and only if $z_u\ge0$ and $z_v\ge0$.

For the arbitrary segment $I=[a,b]$ of the length $d=b-a$ the monotonically
non-decreasing functions are correctly defined for all $x\in I$:
$$
S_u(x)=\esssup\limits_{t\in[a,x]}\bigl(
s_u(t)+g_u(t)
\bigr),\quad S_v(x)=\esssup\limits_{t\in[a,x]}\bigl(
s_v(t)+g_v(t)
\bigr)
$$

We define the function $S(x)=(S_u(x),S_v(x))$, $x\in I$. For arbitrary $y\in\overset{\circ}{W_2^1}(I)$
$$
-\int\limits_I S(x)\,d|y(x)|^2=e^{i\alpha}\int\limits_I |y(x)|^2 dS_u(x)+e^{i\beta}\int\limits_I |y(x)|^2 dS_v(x)\in\Lambda.
$$

Denote $\sigma(x)=S(x)-s(x)$, $x\in I$. It follows from the condition of the theorem:
$$
g_u(x)\le\sigma_u(x)\le f_u(x),\quad
g_v(x)\le\sigma_v(x)\le f_v(x),
$$
thus there exists $C>0$ common to all segments of length $d$, and $\|\sigma\|_I\le C$
uniformly over all segments $|I|=d$ in $L_2(I)$ metric. Note that $\sigma$ itself depends on $I$, while 
the upper estimate of $\|\sigma\|_I$ depends on $d$ only.

Applying the Cauchy--Bunyakovsky inequality,
\begin{equation}
\label{eqsigmaeval}
\begin{aligned}
\left|\int\limits_I \sigma(x)\,d|y(x)|^2\right|=
\left|\int\limits_I \sigma(x)\bigl(y\,\overline{y}'+\overline{y}\,y'\bigr) \,dx\right|
\le\\
\le2\max_{x\in I}|y(x)|\,\|\sigma\|_I\,\|y'\|_I\le2\sqrt{d}\,C\|y'\|_I^2.
\end{aligned}
\end{equation}

Before applying lemma \ref{lmlocalprinc}, we reduce $d$ so that for any $y\in\overset{\circ}{W_2^1}(I)$:
$$
\mathfrak{l}\vint{I}{}[y]=\|y'\|_I^2+\int\limits_I \sigma(x)\,d|y(x)|^2-\int\limits_I S(x)\,d|y(x)|^2\in \Lambda.
$$
For this, we consider two cases:
\begin{itemize}
\item If $\alpha<0<\beta$, we take $0<\varepsilon<\pi/4$, so that $\alpha<-\varepsilon<0<\varepsilon<\beta$.
We decrease $d$ so that $2C\sqrt{d}<\tan\varepsilon$. Taking into account \eqref{eqsigmaeval}
uniformly over all segments of length $d$, the first two terms of the sum for $\mathfrak{l}\vint{I}{}[y]$ lie in the sector
$\Pi_\varepsilon=\{|\arg z|<\varepsilon\}\subset\Lambda$, thus $\mathfrak{l}\vint{I}{}[y]\in\Lambda$.

\item If $\alpha=\beta=0$, we decrease $d$ so that $2C\sqrt{d}<1$.
\end{itemize}
the sectoriality of $L_0$ follows from the lemma \ref{lmlocalprinc}.

Now we prove the compactness criterion. In view of the theorem \ref{thness}, only the proof of sufficiency is of interest.
Consider $d>0$ and the system of intervals $I_k=[a_k,b_k]$, $|I_k|=b_k-a_k=d$, $a_k\to\infty$ as $k\to\infty$.
Representing $s(x)=S_k(x)-\sigma_k(x)$ on $x\in I_k$, as in the first part of the proof, 
for some $C>0$ the estimate $\|\sigma_k\|_{I_k}<C$ holds for all $k\in\NN$. As in the first part,
let $d$ be small enough that $\mathfrak{l}\vint{I_k}{}[y]\in\Lambda$ for all $k\in\NN$ and $y\in\overset{\circ}{W_2^1}(I_k)$, moreover, let
\begin{equation}
\label{eqtaudlimit}
2C\sqrt{d}<\frac{1}{2}\cos\frac{\alpha+\beta}{2}.
\end{equation}

For any
$0<\delta\le d$ and any system of segments $\Delta_k\subset I_k$, $|\Delta_k|=\delta$
$$
\lim\limits_{k\to\infty}\int\limits_{\Delta_k}\int\limits_{\Delta_k}|S_k(\xi)-S_k(\eta)|^2\,d\xi d\eta=+\infty,
$$
and therefore, $|S_k(\xi)-S_k(\eta)|\to+\infty$ uniformly in $\xi,\eta\in I_k$, $|\xi-\eta|=\delta$. Then
uniformly in $\xi,\eta\in I_k$, $\eta<\xi$, $\xi-\eta=\delta$
$$
V_{k}(\xi)-V_k(\eta)\to+\infty,
$$
where
\begin{equation}
\label{eqvk}
V_{k}(x)=\re \Bigl\{S_k(x)\exp(-i\frac{\alpha+\beta}{2})\Bigr\}=\cos\frac{\beta-\alpha}{2}(S_{ku}(x)+S_{kv}(x))
\end{equation}
--- is a non-decreasing functions. Consider
$$
\tau_k(x)=\re \Bigl\{\sigma_k(x)\exp(-i\frac{\alpha+\beta}{2})\Bigr\},
$$
it follows $\|\tau_k\|_{I_k}<C$ for all $k\in\NN$. 

Now we turn to the form
\begin{gather*}
\re\Bigl\{\mathfrak{l}\vint{I_k}{}[y]\exp(-i\frac{\alpha+\beta}{2})\Bigr\}=\cos\frac{\alpha+\beta}{2}
\|y'\|_{I_k}^2+\int\limits_{I_k} |y(x)|^2\,dV_{k}(x)+\\
+\re \Bigl\{\exp(-i\frac{\alpha+\beta}{2})\int\limits_{I_k} \sigma_k(x)\,d|y(x)|^2\Bigr\}=\\
=\Big\{\frac{1}{2}\cos\frac{\alpha+\beta}{2}\|y'\|_{I_k}^2+\int\limits_{I_k} |y(x)|^2\,dV_{k}(x)\Bigr\}+
\Big\{\frac{1}{2}\cos\frac{\alpha+\beta}{2}\|y'\|_{I_k}^2+\int\limits_{I_k} \tau_k(x)\,d|y(x)|^2
\Bigr\}.
\end{gather*}

By \eqref{eqtaudlimit} the expression in the second pair of braces is non-negative. According to lemma \ref{lmSkmonoton} of Appendix 
the expression in the first pair of braces increases without bound as $k\to+\infty$.

Since $\mathfrak{l}\vint{I_k}{}[y]$ is strictly sectorial, the limit \eqref{eqlimlemmaismag} holds.
Applying the lemma \ref{lmlocalprinc}, we conclude that
$L_0$ has extensions with compact resolvents.\qquad$\Box$
\bigskip

{\noindent\bf Proof of theorems \ref{thismagl} and \ref{thbrink}.} We use the same notation as in the proof of theorem \ref{thcriter}.

Consider $d>0$ and the system of intervals $I_k=[a_k,b_k]$, $|I_k|=b_k-a_k=d$, $a_k\to\infty$ as $k\to\infty$.
Let
$d$ be small enough so that $\mathfrak{l}\vint{I_k}{}[y]\in\Lambda$ for any $k\in\NN$ and $y\in\overset{\circ}{W_2^1}(I_k)$.

On each segment $x\in I_k$ we represent $s(x)=S_k(x)-\sigma_k(x)$.

Arguing as in the theorem \ref{thcriter},
the operator $L_0$ has extensions with compact resolvents if and only if for any system of segments $I_k$ of length $|I_k|=d$
and any $0<\delta\le d$ the value
$|S_k(\xi)-S_k(\eta)|\to+\infty$ uniformly in $\xi,\eta\in I_k$, $|\xi-\eta|=\delta$ 
as $k\to\infty$.

The following conditions are equivalent:
\begin{itemize}
\item For any $A>0$
$$
\lim_{k\to+\infty}\mu\set{(\xi,\eta)\in I_k\times I_k}{|S_{k}(\xi)-S_k(\eta)|<A}=0,
$$
\item Uniformly in $\xi,\eta\in I_k$, $\xi-\eta=\delta$, where $0<\delta\le d$
$$
|S_{k}(\xi)-S_k(\eta)|\to+\infty,\mbox{ as }k\to\infty,
$$
\end{itemize}
which follows from the equivalence of the corresponding expressions if we replace the sectorial complex-valued $S_k$ with monotone real-valued $V_k$ \eqref{eqvk}.

In the condition of the theorem \ref{thismagl}, uniformly in $k\in\NN$ holds $\|\sigma_k\|_{I_k}<C$, thus for any $B>0$
\begin{gather*}
\mu\set{x\in I_k}{|\sigma_k(x)|\ge B}\le\frac{C}{B^2},\\
\mu\set{(\xi,\eta)\in I_k\times I_k}{|\sigma_k(\xi)|\ge B\mbox{ or }|\sigma_k(\eta)|\ge B}\le\frac{2dC}{B^2}.
\end{gather*}

Then
\begin{gather*}
\mu\set{(\xi,\eta)\in I_k\times I_k}{|s(\xi)-s(\eta)|<A}\le\\
\le\frac{2dC}{B^2}+\mu\set{(\xi,\eta)\in I_k\times I_k}{|S_{k}(\xi)-S_k(\eta)|<A+2B}\le\\
\le\frac{4dC}{B^2}+\mu\set{(\xi,\eta)\in I_k\times I_k}{|s(\xi)-s(\eta)|<A+4B},
\end{gather*}
which proves theorem \ref{thismagl}.

In the condition of the theorem \ref{thbrink}, on the other hand, uniformly in $k\in\NN$ and $x\in I_k$ holds $\max|\sigma_k(x)|<C$, then
uniformly in $\xi,\eta\in I_k$, $\xi-\eta=\delta$, where $0<\delta\le d$:
$$
|s(\xi)-s(\eta)|\le |S_{k}(\xi)-S_k(\eta)| +2C \le |s(\xi)-s(\eta)| + 4C,
$$
which proves theorem \ref{thbrink}.\qquad$\Box$
\bigskip

{\noindent\bf Proof of theorem \ref{thunifpert}.} Expanding $\Lambda_{\alpha,\beta}(z_0)$ if necessary,
we assume that either
$\alpha=\beta=0$ or $-\pi<\alpha<0<\beta<\pi$.

In accordance with Lemma \ref{lmw2unif} of Appendix, we represent:
$$
\sigma(x)=\sigma_1(x)+\gamma(x),
$$
where for some some $C>0$: $\|\sigma_1\|_I<C$ and for all $x,y\in I$ $|\gamma(x)-\gamma(y)|<C$ uniformly over all intervals $I$, $|I|\le 1$.

For an arbitrary interval $I$, and $y\in \mathscr{D}(\mathfrak{l}_I)=\overset{\circ}{W_2^1}(I)$ we transform the form $\mathfrak{l}_I$ as:
\begin{gather*}
\mathfrak{l}_I[y]=\int\limits_I |y'(x)|^2\,dx - \int\limits_I (s(x)+\sigma(x))\,d|y(x)|^2=\\
=\frac{1}{1+\varepsilon}\left\{
\int\limits_I |y'(x)|^2\,dx - \int\limits_I (1+\varepsilon)s(x)\,d|y(x)|^2
\right\}+\\
+\frac{\varepsilon}{2(1+\varepsilon)}
\left\{
\int\limits_I |y'(x)|^2\,dx - \int\limits_I
\frac{2(1+\varepsilon)}{\varepsilon}\sigma_1(x)\,d|y(x)|^2
\right\}+\\
+\frac{\varepsilon}{2(1+\varepsilon)}
\left\{
\int\limits_I |y'(x)|^2\,dx - \int\limits_I
\frac{2(1+\varepsilon)}{\varepsilon}\gamma(x)\,d|y(x)|^2
\right\}.
\end{gather*}

By the condition of the theorem, the term in the first pair of braces lies in the sector $\Lambda_{\alpha,\beta}(z_0)$. Arguing as in the proof of theorem \ref{thcriter},
we find $0<\delta<1$ so that for all segments $I$ of length $|I|=\delta$ the terms of the second and third pairs of braces:
\begin{itemize}
\item lie in the sector
$\{|\arg z|<\varepsilon_1\}$ for some small $\varepsilon_1>0$ such that $\alpha<-\varepsilon_1<0<\varepsilon_1<\beta$ (in the case $\alpha<0<\beta$);
\item are non-negative (in the case $\alpha=\beta=0$).
\end{itemize}

Applying lemma \ref{lmlocalprinc} to the arbitrary system of segments of length $\delta$, we complete the proof.\qquad$\Box$

\section{General form of the antiderivative of the potential for the semi-bounded operator, counterexamples}
\label{sect5}

{\noindent\bf Proof of lemma \ref{lmmiura}.}
For $y\in W_{2,loc}^1(\RR_+)$ with compact support such that $y(0)=0$ and arbitrary real-valued $\gamma\in L_{2,loc}(\RR_+)$
$$
\int\limits_0^{+\infty}|y'(x)|^2\,dx-\int\limits_0^{+\infty}\Bigl(
\gamma(x)+\int\limits_0^x\gamma^2(t)\,dt
\Bigr)
\,d|y(x)|^2=\int\limits_0^{+\infty}|y'(x)-\gamma(x)y(x)|^2\,dx\ge0,
$$
thus only proof of necessity is of interest.

Let $s$ be a generalized antiderivative of the potential for some semi-bounded operator. Denote by $\hat s$ the extension of $s$ by zero on $(-1,+\infty)$.

The minimal operator $L_0$ in $L_2(-1,+\infty)$,
constructed by $\hat s$ with the differential expression $l$ of the form \eqref{eqly}, will be semi-bounded by virtue of the lemma \ref{lmnormcrit}.

Adding to $\hat s$ the term $-bx$ if needed, we assume that
$L_0>I$.

Now we take the non-trivial solution $\varphi$ of the equation $l(y)=0$, for which $\varphi(-1)=0$. The function $\varphi$ has no zeros for $x>-1$,
otherwise, if there exists $x_0>-1$, $\varphi(x_0)=0$, consider $\psi$, which coincides with $\varphi$ for $x\in[-1,x_0]$ and $\psi(x)=0$ for $x>x_0$,
the minimal quadratic form given by the expression of the form \eqref{eqformmS} $\mathfrak{l}[\psi]=0$, which, taking into account
the statement \ref{stsect}, contradicts the fact that $L_0>I$.

Further for $x\in\RR_+$ let:
$$
\gamma(x)=\hat s(x)+\frac{\varphi^{[1]}(x)}{\varphi(x)}=s(x)-bx+\frac{\varphi^{[1]}(x)}{\varphi(x)},
$$
we find the desired function.\qquad$\Box$
\bigskip

{\noindent\bf Proof of theorem \ref{ths1s2exist}.} We are going to construct both counterexamples in the form:
$$
s_j(x)=S_j(x)+g_j(x)+\int g_j^2(t)\,dt,\quad j=1,2,
$$
where $S_j$ --- monotonically non-decreasing, $S_j'(x)\to+\infty$ as $x\to+\infty$, at the points of differentiability of $S_j$, $g_j\in L_{2,loc}(\RR_+)$.
According to lemmas \ref{lmmiura} and \ref{lmlocalprinc}, the operators $L_{0,j}$ constructed by $s_j$ are semi-bounded
and have extensions with compact resolvents.

Let $h>1+\sqrt{7}$, $\delta=3/h^2$ (in this case $0<\delta<1/2$ and $1/h<1/2-\delta$). We introduce two auxiliary functions. 
$$
v_h(x)=
\begin{cases}
x,&\ x\in[0,1/2-\delta],\\
1/2-\delta+h(x+\delta-1/2)/\delta,&\ x\in [1/2-\delta,1/2],\\
v_h(1-x),&\ x\in[1/2,1],
\end{cases}
$$
$$
\gamma_h(x)=
\begin{cases}
h,&\ x\in[0,1/h],\\
v_h'(x)/v_h(x),&\ x\in[1/h,1/2),\\
-\gamma_h(1-x),&\ x\in(1/2,1].
\end{cases}
$$

Obviously $v_h\in\overset{\circ}{W_2^1}[0,1]$, $\gamma_h\in L_2[0,1]$. For arbitrary $\alpha>0$ the form
\begin{gather*}
\mathfrak{m}[v_h]=\int\limits_0^{1}\bigl(v_h'(x)-\gamma_h(x)v_h(x)\bigr)^2\,dx-
\alpha\int\limits_0^{1}\gamma_h^2(x)v_h^2(x)\,dx\le\\
\le\frac{2}{3h}-2\alpha\int\limits_{1/2-\delta}^{1/2}(v_h'(x))^2\,dx=\frac{2}{3h}-\frac{2\alpha}{3}h^4,
\end{gather*}
at the same time:
$$
\|v_h\|_{I}^2=25/24+O(1/h),\mbox{ as } h\to+\infty,
$$
thus:
$$
\inf\limits_{
\begin{smallmatrix}y\in\overset{\circ}{W_2^1}[0,1]\\
\|y\|_{[0,1]}=1
\end{smallmatrix}
}\mathfrak{m}[y]\to -\infty,\mbox{ as }h\to+\infty.
$$

First we construct $s_2$. To simplify notation, we omit the index $j=2$. Let $\theta\in(0,1)$ be chosen arbitrarily, set $(1+\varepsilon)=1/(1-\theta)^{1/3}$,
$\alpha=1-(1-\theta)^{1/3}$.

We are interested in the properties of the minimal operator $L_{0,\varepsilon}$ constructed by the antiderivative $(1+\varepsilon)s(x)$.
For arbitrary segment $I$, $|I|=1$ and $y\in\overset{\circ}{W_2^1}(I)$ consider:
\begin{gather}
\mathfrak{l}_I[y]= \int\limits_I|y'(x)|^2\,dx-
\int\limits_I(1+\varepsilon)s(x)\,d|y(x)|^2=
\notag
\\
=(1+\varepsilon)\int\limits_I |y(x)|^2S'(x)\,dx+
\int\limits_I\bigl|y'(x)-(1+\varepsilon)g(x)y(x)\bigr|^2\,dx-
\label{eqth5prfform}
\\
-\bigl((1+\varepsilon)^2-(1+\varepsilon)\bigr)
\int\limits_I g^2(x)|y(x)|^2\,dx=(1+\varepsilon)\int\limits_I |y(x)|^2S'(x)\,dx+
\mathfrak{m}_I[y],
\notag
\end{gather}
where
$$
\mathfrak{m}_I[y]=\int\limits_I\bigl|y'(x)-(1+\varepsilon)g(x)y(x)\bigr|^2\,dx-\bigl((1+\varepsilon)^2-(1+\varepsilon)\bigr)
\int\limits_I g^2(x)|y(x)|^2\,dx.
$$

We define $g\equiv 0$ everywhere except for the intervals $[N_k,N_k+1]$, $k\in\NN$, defining
$g(x)=(1-\theta)^{1/3}\gamma_{h_k}(x-N_k)$ for $x\in[N_k,N_k+1]$. The numbers $N_k\in\NN$ and $h_k>0$ will be chosen later.
We denote the segments
$I_{1,k}=[N_k-1/2,N_k+1/2]$, $I_{2,k}=[N_k,N_k+1]$, $I_{3,k}=[N_k+1/2,N_k+3/2]$.

For $y_{h_k}=v_{h_k}(x-N_k)$:
$$
\mathfrak{m}_{I_{2,k}}[y_{h_k}]=\mathfrak{m}[v_h].
$$
Denote
$$
\mu_k=\min\limits_{j=1,2,3}
\inf\limits_{
\begin{smallmatrix}y\in\overset{\circ}{W_2^1}(I_{j,k})\\
\|y\|_{I_{j,k}}=1
\end{smallmatrix}
}\mathfrak{m}_{I_{j,k}}[y].
$$
From the observation above, $\mu_k\to-\infty$ as $h_k\to+\infty$ with $k\to\infty$.

We take the sequences $h_k\to+\infty$, $N_k\in\NN$ so that $A_k=-\mu_k/(1+\varepsilon)$
lie in intervals $N_k-1\le A_k<N_k+1$, and $N_{k}-N_{k-1}\ge3$ for all $k\in\NN$.

Let $S'(x)=l-1$ for any interval $x\in(l-1,l)$, $l\in\NN$ that does not have common points with segments of the form $[N_k-1,N_k+2]$.
For $x\in (N_k-1,N_k+2)$ we define
$$
S'(x)=
\begin{cases}
N_k-1,&\ x\in(N_k-1,N_k-1/2),\\
A_k,&\ x\in(N_k-1/2, N_k+3/2),\\
N_k+1,&\ x\in(N_k+3/2, N_k+2).
\end{cases}
$$

Due to the note to lemma \ref{lmlocalprinc}, 
it is sufficient to verify the conditions of lemma \ref{lmlocalprinc} only on intervals of the form $[l-1,l]$, $[l-1/2,l+1/2]$, $l\in\NN$. 
Let $I$ be the segment of this type.

If $I$ does not coincide with any of $I_{1,k}$, $I_{2,k}$, $I_{3,k}$, $k\in \NN$, then
$\mathfrak{m}_I[y]=\|y'\|_I^2$ and the form $\mathfrak{l}_I$ \eqref{eqth5prfform} is certainly semi-bounded. Otherwise, for $y\in\overset{\circ}{W_2^1}(I)$ and $\|y\|_I=1$:
$$
\mathfrak{l}_I[y]=(1+\varepsilon)A_k+\mathfrak{m}_I[y]\ge(1+\varepsilon)A_k+\mu_k=0,
$$
where $0$ is the greatest lower bound of the form $\mathfrak{l}_I$ by construction.

Applying Lemma \ref{lmlocalprinc} we conclude that $L_{0,\varepsilon}$ is semibounded but does not have extentions with compact resolvents.

To construct $s_1$ it suffices to take $S(x)=x^2/2$, $N_k=4+k$, $k\in\NN$, $h_k=N_k$. The construction of the function $g$ is 
similar to the previous case. The forms $\mathfrak{l}_I$ are not uniformly semibounded on intervals $I_k=[N_k,N_k+1]$, which completes
the proof.\qquad$\Box$
\bigskip

\section{Appendix}
\label{scappendix}

\begin{Lemma}
\label{lmDdenseW21}
Let $s\in L_2[0,a]$, denote $y^{[1]}=y'-sy$ for arbitrary $y\in AC[0,a]$. The lineal
$$
\mathscr{D}=\set{y\in L_2[0,a]}{y, y^{[1]}\in AC[0,a],\ y(0)=y^{[1]}(0)=y(a)=y^{[1]}(a)=0}
$$
is dense in the sense of $W_{2}^1[0,a]$ metric in the space
$$
\overset{\circ}{W_2^1}[0,a]=\set{y\in W_2^1[0,a]}{y(0)=y(a)=0}.
$$
\end{Lemma}

{\noindent\bf Proof.} We perform an auxiliary construction: for arbitrary $\delta>0$ and $A\in\CC$ we construct
$\omega_\delta\in AC[0,a]$ with the properties: $\|\omega_\delta\|_2<\delta$,
$\omega_\delta(0)=A$, $\omega_\delta(a)=0$, $\omega_\delta\perp\psi$,
where
$$
\psi(x)=\exp\int_x^a\overline s(t)\,dt
$$
--- the solution to the homogeneous equation $\psi'+\overline s\psi=0$. The bar in $\overline s$ corresponds to complex conjugation. By $\|\cdot\|_2$ we denote the norm in $L_2[0,a]$.

If $A=0$, we set $\omega_\delta\equiv0$.

Let $A\ne0$, for some small $\varepsilon>0$ we set:
$$
\omega_\delta(x)=
\left\{
\begin{array}{cl}
A+kx&\mbox{ for } x\in [0,\varepsilon/2]\\
(2A/\varepsilon+k)(\varepsilon-x)&\mbox{ for } x\in [\varepsilon/2,\varepsilon]\\
0&\mbox{ for } x\in[\varepsilon,a]
\end{array}
\right.,
$$
selecting $k$ from the orthogonality:
$$
\int\limits_0^{\varepsilon/2}(A+kx)\overline{\psi(x)}\,dx+
\int\limits_{\varepsilon/2}^{\varepsilon}(2A/\varepsilon+k)(\varepsilon-x)\overline{\psi(x)}\,dx=0,
$$
what is equivalent to
$$
k\Bigl(
\int\limits_0^{\varepsilon/2}x\overline{\psi(x)}\,dx+\int\limits_{\varepsilon/2}^{\varepsilon}(\varepsilon-x)\overline{\psi(x)}\,dx
\Bigr)=
-A\Bigl(
\frac{2}{\varepsilon}\int\limits_{\varepsilon/2}^{\varepsilon}(\varepsilon-x)\overline{\psi(x)}\,dx+\int\limits_0^{\varepsilon/2}\overline{\psi(x)}\,dx
\Bigr).
$$

Since $\overline{\psi(0)}= C\ne0$, the last equality, due to continuity of $\psi$, is asymptotically
$$
kC\frac{\varepsilon^2}{4}\sim-\frac{3}{4}AC\varepsilon,\mbox{ or }k\sim\frac{3}{\varepsilon}A\ \mbox{as }\varepsilon\to0+,
$$

Thus:
\begin{gather*}
\|\omega_\delta\|_2^2=\int\limits_0^{\varepsilon/2}|A+kx|^2\,dx+|2A/\varepsilon+k|^2\int\limits_{\varepsilon/2}^{\varepsilon}(\varepsilon-x)^2\,dx=\\
=|A|^2\varepsilon
\int\limits_0^{1/2}(1+3x)^2\,dx\,(1+o(1))+|A|^2\,\frac{25}{24}\varepsilon\,(1+o(1))\sim |A|^2\frac{8}{3}\varepsilon.
\end{gather*}

By choosing a sufficiently small $\varepsilon>0$, we achieve $\|\omega_\delta\|_2<\delta$.

Now we proceed directly to the proof of the lemma.

Let $y_0\in\overset{\circ}{W_2^1}[0,a]$,
then $f_0=y_0'-sy_0\in L_2[0,a]$. All $f\in L_2[0,a]$, which can be represented in the form $f=y'-sy$
for some $y\in\overset{\circ}{W_2^1}[0,a]$, are completely characterized by the orthogonality
to $\psi$ in $L_2[0,a]$ metric.

The lemma will be proved if for any $\varepsilon>0$ we find $f\in AC[0,a]$ orthogonal to $\psi$, so that
$f(0)=f(a)=0$, $\|f-f_0\|_2<\varepsilon$. 
In this case the solution $y\in \mathscr{D}$ to the Cauchy problem $y'-sy=f$ with zero initial condition $y(0)=0$ will be the desired one:
due to the boundedness of the Cauchy operator on the segment, $\|y-y_0\|_2<C_1\|f-f_0\|_2<C_1\varepsilon$, 
$\|y'-y_0'\|_2<\|f-f_0\|_2+\|s(y-y_0)\|_2<\varepsilon +\|s\|_2C_1\varepsilon<C_2\varepsilon$.

To construct the desired $f\in AC[0,a]$ we take an arbitrary absolutely continuous function $g$, approximating $f_0$, $\|g-f_0\|_2<\varepsilon/4$ and
denote
$$
g^\perp=g-\psi\frac{(g,\psi)}{(\psi,\psi)}.
$$

Since $f_0$ is orthogonal to $\psi$, then $\|g^\perp-f_0\|<\varepsilon/2$. Now we construct absolutely continuous functions $\omega_1$, $\omega_2$
orthogonal to $\psi$ with $L_2$-norms less than $\varepsilon/4$, $\omega_1(0)=-g^\perp(0)$, $\omega_1(a)=0$, $\omega_2(a)=-g^\perp(a)$, $\omega_2(0)=0$, finally we take
$f=g^\perp+\omega_1+\omega_2$ as the desired function, $\|f-f_0\|<\varepsilon$.\qquad$\Box$

\begin{Lemma}
\label{lmSkmonoton}
Given the sequence of real-valued non-decreasing functions $S_k$ on the interval $I=[a,b]$ and the quadratic forms for
$y\in\overset{\circ}{W_2^1}(I)$ in $L_2(I)$:
\begin{equation}
\label{eqfrmlkI}
\mathfrak{l}_k[y]=\int\limits_I|y'(x)|^2\,dx + \int\limits_I |y(x)|^2\,dS_k(x).
\end{equation}
The condition
\begin{equation}
\label{eqcondliminflm8}
\lim\limits_{k\to\infty} \inf\limits_{
\begin{smallmatrix}y\in\overset{\circ}{W_2^1}(I)\\
\|y\|_{I}=1
\end{smallmatrix}
}\bigl|\mathfrak{l}_k[y]\bigr|=+\infty.
\end{equation}
is equivalent to the following: for any $0<d\le b-a$ the sequence $S_k(y)-S_k(x)\to+\infty$ uniformly in $y-x=d$, $x,y\in I$ as $k\to\infty$.
\end{Lemma}

{\noindent\bf Proof.} First we prove the direct statement. Assume the contrary: for some $d>0$ and $C>0$ there exist segments
$I_k=[a_k,b_k]\subset I$, $b_k-a_k=d$, $S_k(b_k)-S_k(a_k)<C$.

Without loss of generality, we can assume that $I_k\equiv I$. Otherwise, by shifting the variable $\xi=x-a_k$, we reduce the problem to the one on the segment $I=[0,d]$,
where \eqref{eqcondliminflm8} is also satisfied.

Consider $V_k(x)=S_k(x)-S_k(a)$, by Helly's second theorem there is some subsequence $V_{k_s}(x)\to V(x)$ at each point $x\in I$.
Then for every $y\in\overset{\circ}{W_2^1}(I)$, $\|y\|_I=1$ by Helly's first theorem:
$$
\lim_{s\to\infty}\mathfrak{l}_{k_s}[y]=\int\limits_I|y'(x)|^2\,dx + \int\limits_I |y(x)|^2\,dV(x)<+\infty,
$$
which contradicts \eqref{eqcondliminflm8}.

Now we prove the converse by contradiction. Let $S_k(y)-S_k(x)\to+\infty$ uniformly in $y-x=d$, $x,y\in I$ for any $0<d\le b-a$.
Suppose,
there exist numbers $n_k\in \NN$, $n_k\to\infty$,
the sequence $y_{n_k}\in\overset{\circ}{W_2^1}(I)$,
$\|y_{n_k}\|_{I}=1$, such that $\mathfrak{l}_{n_k}[y_{n_k}]<C$. Without loss of generality, we assume that $n_k=k$.
Since each term in \eqref{eqfrmlkI} is non-negative, then $\|y_k'\|_{I}<C$.

For any $\varepsilon>0$, we consider $E_k=\set{x\in I}{|y_k(x)|>\varepsilon}$. Clearly,
$$
\varepsilon^2\mu_{S_k}(E_k)=\varepsilon^2\int\limits_{E_k} dS_k(x)\le\int\limits_I |y(x)|^2\,dS_k(x)<C,
$$
where $\mu_{S_k}$ denotes the measure on $I$ given by the monotone function $S_k$.

The sets $E_k$ are open and can be represented as a union of pairwise disjoint intervals $E_k=\cup_m \Delta_{k,m}$. From the obtained estimate
and the conditions of the lemma, we obtain that for interval lengths $\max_m|\Delta_{k,m}|\to0$ as $k\to+\infty$.

Let $M_k=|y_k(x_{0,k})|=\max|y_k(x)|$ be the maximum for $x\in I$. According to our observation, for any $\delta>0$, for sufficiently large $k>k_0$, 
one can find $\xi_{0,k}\not\in E_k$ such that $|x_{0,k}-\xi_{0,k}|<\delta$, and then
\begin{gather*}
M_k^2-\varepsilon^2\le|y_k(x_{0,k})|^2-|y_k(\xi_{0,k})|^2\le|y_k^2(x_{0,k})-y_k^2(\xi_{0,k})|\le\\
\le \Bigl|\,\int\limits_{\xi_{0,k}}^{x_{0,k}}2M_k|y_k'(x)|\,dx\, \Bigr|\le
2M_k\sqrt{\delta}\|y_k'\|_I\le 2\sqrt{\delta}\|y_k'\|_I^2\sqrt{b-a} \le 2C^2 \sqrt{\delta}\sqrt{b-a}.
\end{gather*}

Since $\varepsilon>0$ and $\delta>0$ are arbitrary, we conclude that $M_k\to0$, which contradicts the condition $\|y_k\|_{I}=1$.
\qquad$\Box$

\begin{Lemma}
\label{lmw2unif}
The distribution $q\in W_{2,unif}^{-1}(\RR_+)$ if and only if one of the following equivalent statements holds for its generalized antiderivative $s$:
\begin{itemize}
\item
The following estimate holds:
\begin{equation}
\label{eqw2unifconst}
\begin{aligned}
&\sup_{\phantom{_+}x\in\RR_+}\min_{c\in\CC}\int\limits_x^{x+1}|s(\xi)-c|^2\,d\xi<+\infty,\mbox{ which is equivalent to}
\\
&\sup_{\phantom{_+}x\in\RR_+}\left\{\int\limits_x^{x+1}|s(\xi)|^2\,d\xi-\Bigl|
\int\limits_x^{x+1}s(\xi)\,d\xi
\Bigr|^2\right\}<+\infty.
\end{aligned}
\end{equation}
\item
There is a representation:
\begin{equation}
\label{eqrepresunif1}
s(x)=\sigma_1(x)+\gamma(x),
\end{equation}
where $\sigma_1\in L_{2,unif}(\RR_+)$, and $\gamma$ is a piecewise constant function, $\gamma(x)\equiv C_l\in\CC_k$ for $l-1\le x<l$, $l\in\NN$, and $|C_{l+1}-C_l|<C$ for some common
constant $C\ge0$.
\item
There is a representation:
\begin{equation}
\label{eqrepresunif2}
s(x)=\sigma_2(x)+\int\tau(x)\,dx,
\end{equation}
where $\sigma_2\in L_{2,unif}(\RR_+)$, and $\tau\in L_{1,unif}(\RR_+)$.
\end{itemize}
\end{Lemma}
{\noindent\bf Proof.} We show that for any $a>0$, $x\ge0$
\begin{equation}
\label{equnifw2eq}
\begin{split}
a\,\min_{c\in\CC}\int\limits_x^{x+a}|s(\xi)-c|^2\,d\xi&\le
\int\limits_x^{x+a}\int\limits_x^{x+a}|s(\xi)-s(\eta)|^2\,d\xi d\eta\le\\
&\le4a\,\min_{c\in\CC}\int\limits_x^{x+a}|s(\xi)-c|^2\,d\xi.
\end{split}
\end{equation}

Adding and subtracting an arbitrary constant $c\in\CC$, we obtain the estimate:
$$
\Bigl(\int\limits_x^{x+a}\int\limits_x^{x+a}|s(\xi)-s(\eta)|^2\,d\xi d\eta)\Bigr)^{1/2}\le
2\Bigl(a\min_{c\in\CC}\int\limits_x^{x+a}|s(\xi)-c|^2\,d\xi\Bigr)^{1/2},
$$
on the other hand,
$$
\min_{c\in\CC}\int\limits_x^{x+a}|s(\xi)-c|^2\,d\xi\le\int\limits_x^{x+a}|s(\xi)-s(\eta)|^2\,d\xi,
$$
integrating over $\eta\in[x,x+a]$ completes the proof \eqref{equnifw2eq}.

The equivalence of the estimates in \eqref{eqw2unifconst} follows from the explicit form of the constant on which the minimum is achieved:
$$
c_{min}=\int\limits_x^{x+1}s(\xi)\,d\xi.
$$

We turn to representations \eqref{eqrepresunif1} and \eqref{eqrepresunif2}. Clearly, each of \eqref{eqrepresunif1} and \eqref{eqrepresunif2} implies
\eqref{eqw2unifequivMain}. Let us show the converse for $q\in W_{2,unif}^{-1}(\RR_+)$, starting with \eqref{eqrepresunif1}.

For each $l\in\NN$, we define:
$$
\gamma(x)\equiv C_l=\int\limits_{l-1}^l s(\xi)\,d\xi,\quad l-1\le x<l.
$$
For an arbitrary $c\in\CC$
$$
|C_{l+1}-C_l|\le\int\limits_{l}^{l+1}|s(\xi)-c|^2\,d\xi+\int\limits_{l-1}^l|s(\xi)-c|^2\,d\xi,
$$
and the last sum is uniformly bounded in $l\in\NN$.

Consider $\sigma_1=s-\gamma$, we have to show that $\sigma_1\in L_{2,unif}(\RR_+)$. The integral 
$$
\int\limits_{l-1}^l|\sigma_1(\xi)|^2\,d\xi= \int\limits_{l-1}^l|s(\xi)|^2\,d\xi-\Bigl|
\int\limits_{l-1}^l s(\xi)\,d\xi\Bigr|^2
$$ 
is uniformly bounded in $l\in\NN$ due to \eqref{eqw2unifconst}.

Now we prove \eqref{eqrepresunif2}. Consider the function $\beta(x)$. At $x\in\NN$ we define $\beta(x-1)=\gamma(x-1)$ for $x\in\NN$, at the remaining points 
$x\in\RR_+\setminus\NN$ we define $\beta(x)$ as a
continuous broken line. Due to the properties of $\gamma$, it is clear that $\tau=\beta'\in L_{1,unif}(\RR_+)$.

Since $|\gamma(x)-\beta(x)|$ is uniformly bounded, $\gamma-\beta\in L_{2,unif}(\RR_+)$. Finally, we define $\sigma_2=\sigma_1+\gamma-\beta$ and
complete the proof.\qquad$\Box$

\end{document}